\documentclass[11pt]{amsart}
\usepackage[T1]{fontenc}
\usepackage{lmodern}
\usepackage{amsmath,amssymb,amsthm}
\usepackage[a4paper,margin=27mm]{geometry}
\usepackage[colorlinks=true,linkcolor=blue,citecolor=blue,urlcolor=blue]{hyperref}
\hypersetup{pdftitle={Polynomial approximation in Fock space},
  pdfauthor={Shengzhao Hou and Liankuo Zhao}}

\newcommand{\C}{\mathbb C}
\newcommand{\F}{\mathcal F}
\newcommand{\M}{\mathcal M}
\newcommand{\R}{\mathbb R}
\newcommand{\Z}{\mathbb Z}

\newcommand{\dA}{\,dA}
\DeclareMathOperator{\Hol}{Hol}
\DeclareMathOperator{\dist}{dist}

\newtheorem{theorem}{Theorem}[section]
\newtheorem{lemma}[theorem]{Lemma}
\newtheorem{proposition}[theorem]{Proposition}
\theoremstyle{remark}
\newtheorem{remark}[theorem]{Remark}
\numberwithin{equation}{section}
\allowdisplaybreaks[2]
\title{Polynomial approximation in Fock space}
\author{Shengzhao Hou}
\address{School of Mathematical Science, Soochow University,
Suzhou 215006, Jiangsu, China}
\email{shou@suda.edu.cn}
\author{Liankuo Zhao}
\address{School of Mathematical Science, Shanxi Normal University,
Taiyuan 030031, Shanxi, China}
\email{lkzhao@sxnu.edu.cn}
\date{}
\makeatletter\@namedef{subjclassname@2020}{2020 Mathematics Subject Classification}\makeatother

\begin{document}
\begin{abstract}
Let $\F_\alpha^2$ be the Fock space with Gaussian weight $e^{-\alpha|z|^2}$,
where $\alpha>0$.
If a nonzero entire function $F$ has quadratic type strictly below $\alpha/2$,
then its polynomial multiples are dense in $F\Hol(\C)\cap\F_\alpha^2$.
This determines the sharp type threshold for a general positive answer
to the Newman--Shapiro approximation problem, without a regularity
hypothesis on the zeros of $F$. The proof uses a
small-value theorem based on a subharmonic comparison principle, followed
by an argument with finitely many translates. At the critical type we
construct an explicit entire function $F_*$ for which the Taylor polynomials
of each admissible entire multiplier give approximation in the weighted norm.
Together with the counterexample of Belov and Borichev, this shows that both
outcomes occur at the critical type.
\end{abstract}
\maketitle

\section{Introduction}

For $\alpha>0$, let $\F_\alpha^2$ denote the Fock space of entire functions $f$
on $\C$ for which
\[
 \|f\|_\alpha^2=\frac{\alpha}{\pi}\int_{\C}|f(z)|^2e^{-\alpha|z|^2}\dA(z)<\infty,
\]
where $dA$ is Lebesgue area measure. We refer to \cite{Zhu} for background
on these spaces.

Let $f$ be a nonzero entire function such that $pf\in\F_\alpha^2$ for every
polynomial $p$. Write $\Hol(\C)$ for the space of entire functions and set
\[
 [f]_\alpha=\overline{\{pf:p\in\C[z]\}}^{\,\F_\alpha^2},
 \qquad \M_\alpha(f)=f\Hol(\C)\cap\F_\alpha^2.
\]
A function in $\F_\alpha^2$ belongs to $\M_\alpha(f)$ exactly when it vanishes
at every zero of $f$ with at least the same multiplicity. Since derivative
evaluations are continuous on $\F_\alpha^2$, $\M_\alpha(f)$ is closed and
$[f]_\alpha\subset\M_\alpha(f)$. The approximation question is whether
\begin{equation}\label{eq-approxquestion}
 [f]_\alpha=\M_\alpha(f).
\end{equation}
Equivalently, one asks whether the polynomials are dense among the entire
functions $V$ with finite norm
\[
 \|V\|_{f,\alpha}^2=\|fV\|_\alpha^2
 =\frac{\alpha}{\pi}\int_{\C}|V(z)|^2|f(z)|^2e^{-\alpha|z|^2}\dA(z).
\]

Newman and Shapiro \cite{NS} posed the problem in their study of adjoints
of multiplication operators, using the exponentials $e^{\lambda z}$ as
test functions. If $e^{\lambda z}f\in\F_\alpha^2$ for every $\lambda\in\C$,
then for each $A>0$ the Fock evaluation estimate \eqref{eq-eval} applied
to $\lambda\in\{\pm\sqrt2 A,\pm i\sqrt2 A\}$ gives
$|f(z)|\le C_A e^{\alpha|z|^2/2-A|z|}$. Under this bound the closed spans
of the exponential and polynomial multiples of $f$ agree
\cite[Lemma~4.1]{BB}, after rescaling to $\alpha=\pi$. The problem is
related to spectral synthesis for systems of reproducing kernels and
to Gaussian Gabor expansions \cite{BBB}. Borichev gave a counterexample
in the Gaussian space with all polynomial multiples in the space, and
obtained counterexamples admitting all exponential multiples for
suitable faster-decaying radial weights \cite[Corollaries~1 and~2]{Borichev}.
Belov and Borichev gave a Gaussian counterexample with all exponential
multiples in the space \cite[Theorem~1.2]{BB} and established positive
results under growth and zero-distribution assumptions.

For a nonzero entire function $f$, its quadratic type is
\[
 \tau_2(f)=\limsup_{r\to\infty}\frac{\log M_f(r)}{r^2},
 \qquad M_f(r)=\max_{|z|=r}|f(z)|.
\]
We set $\tau_2(0)=-\infty$. If $\tau_2(f)<\beta$, then
$|f(z)|\le C_\beta e^{\beta|z|^2}$ for some constant $C_\beta$ and all
$z\in\C$. In particular, if $\tau_2(f)<\alpha/2$, then every polynomial
and exponential multiple of $f$ belongs to $\F_\alpha^2$. The evaluation
estimate \eqref{eq-eval} gives $\tau_2(f)\le\alpha/2$ for every nonzero
$f\in\F_\alpha^2$.

Belov and Borichev proved that \eqref{eq-approxquestion} holds in
$\F_\pi^2$ for functions of sufficiently small quadratic type
\cite[Theorem~1.5]{BB}. The following theorem extends their conclusion to
the full range below the critical type.

\begin{theorem}\label{thm-main}
If $F$ is a nonzero entire function in $\F_\pi^2$ and $\tau_2(F)<\pi/2$, then
\begin{equation}\label{eq-main}
 [F]_\pi=\M_\pi(F).
\end{equation}
In particular, if $V$ is entire and $FV\in\F_\pi^2$, then there are
polynomials $p_n$ such that
\[
 \lim_{n\to\infty}\int_{\C}|F(z)|^2|p_n(z)-V(z)|^2e^{-\pi|z|^2}\dA(z)=0.
\]
\end{theorem}

The change of variables
\begin{equation}\label{eq-scaling}
 F(z)=f\bigl(\sqrt{\pi/\alpha}\,z\bigr)
\end{equation}
is an isometric isomorphism from $\F_\alpha^2$ onto $\F_\pi^2$. It
preserves polynomial and entire multiples and gives
$\tau_2(F)=(\pi/\alpha)\tau_2(f)$. Thus Theorem~\ref{thm-main} is equivalent
to \eqref{eq-approxquestion} for $\tau_2(f)<\alpha/2$ in $\F_\alpha^2$.
The theorem and the Fock evaluation estimate also show that the symbol in
the counterexample of \cite[Theorem~1.2]{BB} necessarily has quadratic type
$\pi/2$ in the normalization of $\F_\pi^2$.

At the critical type, the following product gives a positive example.
It also yields approximation by the Taylor polynomials of the multiplier.

\begin{theorem}\label{thm-endpoint}
The product
\begin{equation}\label{eq-product}
 F_*(z)=\prod_{n=1}^{\infty}\left(1-\left(\frac{z}{n+\sqrt n}\right)^{2n}\right)
\end{equation}
defines an entire function of order two and type exactly $1/2$, and
$e^{\lambda z}F_*\in\F_1^2$ for every $\lambda\in\C$. If $V$ is entire
and $F_*V\in\F_1^2$, then
\begin{equation}\label{eq-taylor}
 \lim_{N\to\infty}\|F_*(V-S_NV)\|_1=0,
\end{equation}
where $S_NV(z)=\sum_{k=0}^N V^{(k)}(0)z^k/k!$ is the $N$th Taylor polynomial
of $V$ at the origin. Consequently, $[F_*]_1=\M_1(F_*)$.
\end{theorem}

To prove Theorem~\ref{thm-main}, suppose that a nonzero multiple
\(g=fV\) is orthogonal to all polynomial multiples of \(f\).
In Section~\ref{sec-main-proof} we associate entire functions \(U\) and
\(S\) with \(g\), with \(\tau_2(U)\le\tau_2(f)\).
Section~\ref{sec-comparison} proves that \(S\) is a polynomial by
combining a small-value theorem with a suitable translate of the
data. This polynomiality yields zeros of \(UV\) near lattice
points. Counting these zeros for finitely many translates gives a
lower bound on quadratic type that is incompatible with
\(\tau_2(f)<\pi/2\). The proof of Theorem~\ref{thm-endpoint} is given
in Section~\ref{sec-critical}. There we compare the multiplier norm
with a radial weighted norm, in which the monomials are orthogonal.

Throughout, $X\asymp Y$ means $C_1Y\le X\le C_2Y$ for positive constants
$C_1,C_2$ independent of the variables under consideration. The letter $C$
denotes a positive constant whose value may change from one occurrence
to the next.

\section{Preliminaries}\label{s2-preliminaries}

Write \(D(z,r)\) for the open disk of center \(z\) and radius \(r\).
For a nonzero entire function \(f\), denote its zero set by \(Z(f)\),
the multiplicity of a zero \(w\) by \(m_f(w)\), and the number of zeros
in \(D(0,r)\), counted with multiplicity, by \(n_f(r)\).

\subsection{Fock-space estimates and quadratic type}

The reproducing kernel and normalized reproducing kernel of
\(\F_\alpha^2\) at \(w\) are
\[
 K_w(z)=e^{\alpha z\overline w},\qquad
 k_w(z)=e^{\alpha z\overline w-\alpha|w|^2/2},
\]
respectively. Since \(f(w)=\langle f,K_w\rangle_\alpha\) and
\(\|K_w\|_\alpha^2=K_w(w)=e^{\alpha|w|^2}\), Cauchy--Schwarz gives
\begin{equation}\label{eq-eval}
 |f(z)|\le\|f\|_\alpha e^{\alpha|z|^2/2}.
\end{equation}
We will use the local estimate
\begin{equation}\label{eq-localeval}
 |f(z)|^2e^{-\alpha|z|^2}
 \le C_\alpha\int_{D(z,1)}|f(w)|^2e^{-\alpha|w|^2}\dA(w).
\end{equation}
The right-hand side tends to zero as \(|z|\to\infty\), since the
integrand is integrable on \(\C\). In particular,
\(|f(z)|e^{-\alpha|z|^2/2}\to0\) at infinity.
We also use the fact that the Taylor polynomials of \(f\) converge
to \(f\) in \(\F_\alpha^2\). Indeed, the normalized monomials
\((\alpha^k/k!)^{1/2}z^k\), \(k\ge0\), form an orthonormal basis, and
the partial sums of the resulting orthogonal expansion are the Taylor
polynomials. See
\cite[Propositions~2.1 and~2.2 and Lemma~2.32]{Zhu}.

Throughout Sections~\ref{s2-preliminaries}--\ref{sec-comparison}, we work
with the normalization \(\alpha=\pi\).

\begin{lemma}\label{lem-product}
If \(f\) and \(g\) are nonzero entire functions with
\(\tau_2(f)+\tau_2(g)<\pi/2\), then \(fg\in[f]_\pi\).
\end{lemma}
\begin{proof}
Choose \(\beta>\tau_2(f)\) and \(\eta>\tau_2(g)\) with
\(\beta+\eta<\pi/2\). The bound \(|f(z)|\le C_\beta e^{\beta|z|^2}\)
ensures that every polynomial multiple of \(f\) belongs to \(\F_\pi^2\).
Similarly, \(|g(z)|\le C_\eta e^{\eta|z|^2}\) and
\(\pi-2\beta-2\eta>0\) show that \(g\in\F_{\pi-2\beta}^2\).
Let \(p_n\) be the \(n\)-th Taylor polynomial of \(g\). Then
\begin{align*}
 \|f(g-p_n)\|_\pi^2
 &=\int_\C |f(z)|^2|g(z)-p_n(z)|^2e^{-\pi|z|^2}\dA(z)\\
 &\le C_\beta^2\int_\C |g(z)-p_n(z)|^2
                  e^{-(\pi-2\beta)|z|^2}\dA(z)\\
 &=\frac{\pi C_\beta^2}{\pi-2\beta}\|g-p_n\|_{\pi-2\beta}^2.
\end{align*}
The last quantity tends to zero by Taylor convergence in
\(\F_{\pi-2\beta}^2\).
\end{proof}

\begin{lemma}[{\cite[Lemma~5.1]{BB}}]\label{lem-Vtype}
There is an absolute constant \(K>1\) such that, if \(f\) and \(V\)
are nonzero entire functions, \(\tau_2(f)<\beta<\pi/2\), and
\(fV\in\F_\pi^2\), then
\[
 \tau_2(V)\le\frac{\pi}{2}+K\beta.
\]
\end{lemma}

\subsection{Thin sets and the square lattice}

Let \(\Lambda=\Z+i\Z\), \(\Lambda_0=\Lambda\setminus\{0\}\), and write
\(N_E(r)=\#(E\cap D(0,r))\) for a locally finite set \(E\subset\C\).
We say that \(E\) has \emph{zero radial density} if \(N_E(r)=o(r^2)\).

Following \cite[Sections~2 and~3]{BB}, we call a measurable set
\(\Omega\subset\C\) \emph{thin} if it can be written as
\(\Omega=\Omega_1\cup\Omega_2\), where
\begin{equation}\label{eq-thin}
 \lim_{r\to\infty}\frac{|\Omega_1\cap D(0,r)|}{r^2}=0,
 \qquad
 \int_{\Omega_2}\frac{\dA(z)}{(1+|z|)^2\log(2+|z|)}<\infty.
\end{equation}
Here \(|\cdot|\) denotes area when applied to a measurable set. Finite
unions and translates of thin sets are thin. For translation, use
\(D(0,r)-a\subset D(0,r+|a|)\) for the first condition and the
comparability of the weight at \(z\) and \(z+a\) for the second, for
each fixed \(a\in\C\).

A thin set satisfies
\begin{equation}\label{eq-lowerdensity}
 \liminf_{r\to\infty}\frac{|\Omega\cap D(0,r)|}{r^2}=0.
\end{equation}
To see this, put \(A_2(t)=|\Omega_2\cap D(0,t)|\) and
\(\omega(t)=(1+t)^{-2}/\log(2+t)\). Because \(\omega(t)\to0\) as
\(t\to\infty\), Tonelli's theorem gives
\[
 \int_{\Omega_2}\omega(|z|)\dA(z)
 =\int_0^\infty A_2(t)\bigl(-\omega'(t)\bigr)\,dt.
\]
If \(\liminf_{t\to\infty}A_2(t)/t^2>0\), then \(A_2(t)\ge ct^2\)
for large \(t\). Since
\(-\omega'(t)\ge 2(1+t)^{-3}/\log(2+t)\), the integral would diverge.
Thus \(A_2(t)/t^2\) has lower limit zero, and the contribution of
\(\Omega_1\) is \(o(t^2)\) on all radii. This proves
\eqref{eq-lowerdensity}.

A set \(\Gamma\subset\Lambda_0\) is \emph{lattice thin} if
\(\bigcup_{\lambda\in\Gamma}D(\lambda,r)\) is thin for every fixed
\(r>0\). The next lemma gives one sequence of radii for a finite
family of translated exceptional sets.

\begin{lemma}\label{lem-commonradii}
Let \(a_0,\ldots,a_{l-1}\in\C\) and let
\(\Gamma_0,\ldots,\Gamma_{l-1}\subset\Lambda_0\) be lattice thin.
Set
\[
 \mathcal E=\bigcup_{j=0}^{l-1}\bigl((\Gamma_j\cup\{0\})-a_j\bigr).
\]
Then there are radii \(r_k\to\infty\) for which
\(N_{\mathcal E}(r_k)=o(r_k^2)\).
\end{lemma}
\begin{proof}
Fix \(0<d<1/2\). The set
\(\Omega=\bigcup_{\lambda\in\mathcal E}D(\lambda,d)\) is thin by
lattice thinness, translation invariance, and closure under finite
unions. The disks at the finitely many translated origins contribute
only bounded sets. Hence
\eqref{eq-lowerdensity} supplies \(s_k\to\infty\) such that
\(|\Omega\cap D(0,s_k)|=o(s_k^2)\). Let \(r_k=s_k-d\).
The disks centered at points of \(\mathcal E\cap D(0,r_k)\) lie in
\(\Omega\cap D(0,s_k)\). They have overlap multiplicity at most \(l\)
because disks within any one of the translated lattices are disjoint.
Integrating the sum of disk indicators gives
\[
 \pi d^2 N_{\mathcal E}(r_k)
 =\int_{\Omega\cap D(0,s_k)}
     \sum_{\lambda\in\mathcal E\cap D(0,r_k)}\mathbf1_{D(\lambda,d)}\,dA
 \le l|\Omega\cap D(0,s_k)|.
\]
The right-hand side is \(o(r_k^2)\) because \(s_k/r_k\to1\).
\end{proof}

Let \(\sigma\) be the Weierstrass sigma function of the square lattice
\(\Lambda\). Its zeros are simple, \(Z(\sigma)=\Lambda\), and
\begin{equation}\label{eq-sigma}
 |\sigma(z)|\asymp\dist(z,\Lambda)e^{\pi|z|^2/2}.
\end{equation}
It has completely regular growth with constant indicator \(\pi/2\).
See \cite[Section~2]{BB}. Set \(\sigma_0(z)=\sigma(z)/z\), with the value
at \(0\) given by the removable extension. For \(|z|\ge1\), dividing
\eqref{eq-sigma} by \(|z|\) and using
\(\dist(z,\Lambda)\le 1/\sqrt2\) gives the upper bound below.
On \(|z|\le1\), it follows by continuity. Thus
\begin{equation}\label{eq-sigma0}
 |\sigma_0(z)|e^{-\pi|z|^2/2}\le\frac{C}{1+|z|},
\end{equation}
and, for fixed \(0<d<1/4\),
\begin{equation}\label{eq-sigma0-reverse}
 |\sigma_0(z)|e^{-\pi|z|^2/2}\ge\frac{c_d}{1+|z|},
 \qquad z\notin\bigcup_{\lambda\in\Lambda} D(\lambda,d).
\end{equation}
For this lower bound the distance to \(\Lambda\) is at least \(d\),
and \(|z|\le 1+|z|\).
Also, uniformly for \(\lambda\in\Lambda_0\) and \(0<s<1/4\),
\begin{equation}\label{eq-circle}
 |\sigma_0(z)|e^{-\pi|z|^2/2}\asymp\frac{s}{1+|\lambda|},
 \qquad |z-\lambda|=s.
\end{equation}
Indeed, on the indicated circle the nearest lattice point is
\(\lambda\), and \(|z|\asymp 1+|\lambda|\).

\subsection{Cauchy transforms}

For \(f,g\in\F_\pi^2\), define the Cauchy transforms
\begin{equation}\label{eq:CAB}
\begin{aligned}
 C_{f,g}(z)&=\int_\C\frac{f(w)\overline{g(w)}}{z-w}
               e^{-\pi|w|^2}\dA(w),\\
 A_g(z)&=C_{g,g}(z),\\
 B_g(z)&=\int_\C\frac{\sigma_0(w)\overline{g(w)}}{z-w}
               e^{-\pi|w|^2}\dA(w).
\end{aligned}
\end{equation}

\begin{lemma}\label{lem-cauchy}
The integrals in \eqref{eq:CAB} converge absolutely and define continuous
functions. The transforms \(C_{f,g}\) and \(A_g\) are bounded and tend
to zero at infinity. Furthermore,
\begin{equation}\label{eq-D}
 |B_g(z)|\le C_g\frac{\sqrt{\log(2+|z|)}}{1+|z|}.
\end{equation}
\end{lemma}
\begin{proof}
These are the estimates of \cite[Lemma~3.1]{BB}. We include the
bounds and continuity argument for the transforms in \eqref{eq:CAB}.
Write \(p_h(w)=|h(w)|e^{-\pi|w|^2/2}\) for \(h\in\F_\pi^2\). Then
\(p_h\in L^2(\C,dA)\), \(\|p_h\|_2=\|h\|_\pi\), and \(p_h\) is
bounded and tends to zero at infinity. The absolute value of the
density of \(C_{f,g}\) is \(p_fp_g\), which is integrable and tends to
zero. On \(|z-w|<1\) the absolute integral is bounded by
\[
 2\pi\sup_{w\in D(z,1)}p_f(w)p_g(w),
 \qquad\text{since}\qquad
 \int_{|\zeta|<1}\frac{\dA(\zeta)}{|\zeta|}=2\pi.
\]
This bound tends to zero as \(|z|\to\infty\). On the complement
\(|z-w|\ge1\), the integral is bounded by \(\|p_fp_g\|_1\). It
also tends to zero by dominated convergence because
\(\mathbf1_{\{|z-w|\ge1\}}/|z-w|\to0\) for each fixed \(w\).
This proves the assertion for \(C_{f,g}\) and, by taking \(f=g\),
for \(A_g\).

For \(B_g\), \eqref{eq-sigma0} bounds the absolute density by
\(C p_g(w)/(1+|w|)\). The part with \(|z-w|<1\) is bounded by
\(C\|g\|_\pi/(1+|z|)\), using \eqref{eq-eval}. On the part
with \(|z-w|\ge1\), Cauchy--Schwarz gives the bound
\[
 C\|g\|_\pi\left(
      \int_{|z-w|\ge1}\frac{\dA(w)}{(1+|w|)^2|z-w|^2}
      \right)^{1/2}.
\]
For \(R=|z|\ge4\), split this last integral into three regions.
On \(|w|\le R/2\), use \(|z-w|\ge R/2\) and integrate
\((1+|w|)^{-2}\) over the disk to get \(O(R^{-2}\log R)\).
On \(|w|>R/2\) and \(1\le |z-w|\le 2R\), use
\((1+|w|)^{-2}\le 4/R^2\) and integrate \(|z-w|^{-2}\)
over the annulus, again obtaining \(O(R^{-2}\log R)\).
On \(|z-w|>2R\), use \(|w|>|z-w|/2\) and integrate
\(4|z-w|^{-4}\) to get \(O(R^{-2})\). For bounded \(z\), the
same integral is bounded uniformly, by splitting at a fixed large
\(|w|\). This proves \eqref{eq-D}.

All these estimates give absolute convergence. For continuity at a
fixed \(z_0\), first restrict to \(w\in D(z_0,2d)\) with
\(|z-z_0|<d\). Local boundedness of the densities and integration of
\(|z-w|^{-1}\) over \(D(z,3d)\) give a bound \(O(d)\) for this part.
Outside \(D(z_0,2d)\), the kernels converge pointwise as \(z\to z_0\)
and are bounded in modulus by \(2/|z_0-w|\). This gives an
integrable majorant by absolute convergence at \(z_0\).
Letting \(d\downarrow0\) proves continuity.
\end{proof}

For \(0\ne g\in\F_\pi^2\), put
\begin{equation}\label{eq-I}
 I_g=A_g-\frac{gB_g}{\sigma_0}.
\end{equation}
This is the diagonal Cauchy transform in \cite[Section~3]{BB}.
In particular,
\(\sigma_0I_g\) extends to an entire function.

\begin{lemma}\label{lem-BB}
Let \(0\ne g\in\F_\pi^2\).
\begin{enumerate}
\item[\rm(i)] There is a thin set \(\Omega\) such that
\begin{equation}\label{eq-Iasymptotic}
 I_g(z)=\frac{\|g\|_\pi^2+o(1)}{z},
 \qquad |z|\to\infty,\quad z\notin\Omega.
\end{equation}
\item[\rm(ii)] If \(I_g=H/\sigma_0+P\), where \(H\) is entire and
\(P\) is a polynomial, then for each \(r>0\) there is a lattice thin
set \(\Gamma\subset\Lambda_0\) such that \(H\) has a zero in every
disk \(D(\lambda,r)\) with \(\lambda\in\Lambda_0\setminus\Gamma\).
\end{enumerate}
\end{lemma}
These are the diagonal cases of \cite[Lemmas~3.2 and~3.5]{BB}.

\begin{lemma}\label{lem-zeros}
For \(0\ne g\in\F_\pi^2\) and any polynomial \(P\), the entire function
\(\sigma_0(I_g+P)\) is nonzero.
\end{lemma}
\begin{proof}
If \(\sigma_0(I_g+P)\equiv0\), then \(I_g=-P\) off \(\Lambda_0\).
Let \(\Omega\) be the thin set from Lemma~\ref{lem-BB}(i).
The complement of \(\Omega\cup\Lambda_0\) is unbounded by
\eqref{eq-lowerdensity}, since \(\Lambda_0\) has area zero. Along a
sequence in this complement with modulus tending to infinity,
\eqref{eq-Iasymptotic} gives
\(zP(z)\to-\|g\|_\pi^2\). This is impossible. For \(P\equiv0\) the
left-hand side is zero, and for \(P\not\equiv0\) its modulus tends to
infinity as \(|z|\to\infty\).
\end{proof}

\section{Proof of the main theorem}\label{sec-main-proof}

\subsection{The auxiliary entire functions}

Suppose the conclusion of Theorem~\ref{thm-main} fails for \(f\).
The spaces \([f]_\pi\subset \M_\pi(f)\) are closed in \(\F_\pi^2\),
so the orthogonal complement of \([f]_\pi\) in \(\M_\pi(f)\) contains
a nonzero vector \(g\). Since \(g\in\M_\pi(f)\), we can write
\(g=fV\) with \(V\) entire. Thus
\begin{equation}\label{eq-defect}
 f\not\equiv0,\qquad \tau_2(f)<\frac\pi2,\qquad
 0\ne g=fV\in\F_\pi^2,\qquad g\perp[f]_\pi.
\end{equation}
Choosing any \(\beta_0\) with
\(\tau_2(f)<\beta_0<\pi/2\), we obtain
\(\tau_2(V)\le\pi/2+K\beta_0\) from Lemma~\ref{lem-Vtype}.
Define
\begin{equation}\label{eq-US}
 U=\sigma_0 C_{f,g}-fB_g,\qquad S=VC_{f,g}-A_g.
\end{equation}
The product \(UV\) will be used to count zeros.
The estimate \(\tau_2(U)\le\tau_2(f)\) proved below bounds the cost
in quadratic type when copies of \(U\) are multiplied together.
Proposition~\ref{prop-polynomial} will show that \(S\) is a
polynomial, allowing us to apply Lemma~\ref{lem-BB}(ii) to \(UV\).

\begin{proposition}\label{prop-entire}
Under \eqref{eq-defect}, \(U\) and \(S\) are entire and
\begin{equation}\label{eq-identity}
 UV=\sigma_0(I_g+S),\qquad
 U\not\equiv0,\qquad\tau_2(U)\le\tau_2(f).
\end{equation}
The function \(S\) is of finite quadratic type.
\end{proposition}
\begin{proof}
By \(\bar\partial(1/z)=\pi\delta_0\), convolution with the densities
in \eqref{eq:CAB} gives the distributional identities
\begin{align*}
 \bar\partial C_{f,g}&=\pi f\overline g\,e^{-\pi|z|^2},\\
 \bar\partial A_g&=\pi|g|^2e^{-\pi|z|^2},\\
 \bar\partial B_g&=\pi\sigma_0\overline g\,e^{-\pi|z|^2}.
\end{align*}
The absolute integral bounds in Lemma~\ref{lem-cauchy} justify
integrating the kernels against compactly supported test
functions. Since \(\sigma_0\), \(f\), and \(V\) are holomorphic,
\begin{align*}
 \bar\partial U
 &=\pi(\sigma_0 f\overline g-f\sigma_0\overline g)e^{-\pi|z|^2}=0,\\
 \bar\partial S
 &=\pi(Vf\overline g-|g|^2)e^{-\pi|z|^2}=0.
\end{align*}
The second equation uses \(Vf=g\). By Weyl's lemma and continuity
of the transforms, \(U\) and \(S\) are entire. Off the zeros of
\(\sigma_0\), we have
\begin{align*}
 UV
 &=\sigma_0 V C_{f,g}-gB_g\\
 &=\sigma_0\left(VC_{f,g}-A_g+A_g-\frac{gB_g}{\sigma_0}\right)
  =\sigma_0(S+I_g).
\end{align*}
Both sides extend across \(\Lambda_0\), giving the identity on
\(\C\). Also, \(C_{f,g}\) and \(A_g\) are bounded, so
\(|S(z)|\le C(1+|V(z)|)\). Since \(\tau_2(V)<\infty\), \(S\)
has finite quadratic type.

Suppose \(U\equiv0\). The product identity gives \(S=-I_g\) off
\(\Lambda_0\). Fix \(0<d<1/4\). At points with
\(\dist(z,\Lambda)\ge d\), \eqref{eq-eval} and
\eqref{eq-sigma0-reverse} give
\[
 \frac{|g(z)|}{|\sigma_0(z)|}
 \le \frac{\|g\|_\pi}{c_d}(1+|z|).
\]
Combining this with \eqref{eq-D} and the boundedness of \(A_g\) gives
\[
 |S(z)|\le C_{g,d}\bigl(1+\sqrt{\log(2+|z|)}\bigr)
 \qquad \dist(z,\Lambda)\ge d.
\]
If \(z\in D(\lambda,d)\) and \(w\in\partial D(\lambda,d)\), then
\(|w|\le |z|+2d\). The maximum modulus principle on each disk
therefore gives the bound on the whole plane after adjusting the
constant. Cauchy's inequalities yield, for each \(m\ge1\),
\[
 \frac{|S^{(m)}(0)|}{m!}
 \le \frac{C_{g,d}\bigl(1+\sqrt{\log(2+r)}\bigr)}{r^m}
 \longrightarrow0\qquad (r\to\infty).
\]
Thus \(S\) is constant. Lemma~\ref{lem-zeros} now gives
\(\sigma_0(I_g+S)\not\equiv0\), contradicting \(UV\equiv0\).

For the auxiliary function used here, the type estimate in
\cite[Lemma~3.4]{BB} follows directly from the Cauchy transforms by
an argument as in \cite[Lemma~4.3]{BB}. Write
\(\beta=\tau_2(f)\), choose \(2\beta/\pi<\varepsilon<1\), and let
\[
 H(z)=\sigma\bigl(\sqrt{1-\varepsilon}\,z\bigr),\qquad
 \eta=\frac{\pi(1-\varepsilon)}{2}.
\]
Its zeros are simple. The upper estimate in \eqref{eq-sigma}
gives \(\tau_2(H)\le\eta\), and the lower estimate along an unbounded
sequence at fixed positive distance from \(Z(H)\) gives equality.
For \(\lambda\in Z(H)\), the quotient
\(H/(\cdot-\lambda)\) is entire and has the same quadratic type as \(H\).
Indeed, \(|z-\lambda|\asymp r\) uniformly on \(|z|=r\) for large \(r\).
Since \(\beta+\eta<\pi/2\), Lemma~\ref{lem-product} gives
\[
 fH\in[f]_\pi,\qquad
 \frac{fH}{\cdot-\lambda}\in[f]_\pi
 \quad(\lambda\in Z(H)).
\]
In particular, \(fH\) and \(g\) are in \(\F_\pi^2\). The
\(\bar\partial\)-derivative of \(HC_{f,g}-C_{fH,g}\) is zero because
both terms have distributional derivative
\(\pi Hf\overline g\,e^{-\pi|z|^2}\). This entire function vanishes
at every \(\lambda\in Z(H)\), since
\[
 C_{fH,g}(\lambda)=-\left\langle\frac{fH}{\cdot-\lambda},g\right\rangle_\pi=0
\]
by orthogonality. Consequently \(T=C_{f,g}-C_{fH,g}/H\) extends to an
entire function.

The zeros of \(H\) form the lattice
\((1-\varepsilon)^{-1/2}\Lambda\). Take pairwise disjoint disks of
the same small radius \(d_H>0\) about these zeros. Outside them,
\eqref{eq-sigma} gives \(|H(z)|\ge c e^{\eta|z|^2}\). There
\[
 |T(z)|\le |C_{f,g}(z)|+C e^{-\eta|z|^2},
\]
which tends to zero at infinity. The maximum modulus principle on
the disks shows that the same decay holds on the whole plane.
Liouville's theorem gives \(T\equiv0\). Substituting
\(C_{f,g}=C_{fH,g}/H\) into \eqref{eq-US} outside the disks gives
\[
 U=\sigma_0\frac{C_{fH,g}}{H}-fB_g.
\]
Choose \(\beta'\) with \(\beta<\beta'<\pi/2-\eta\).
By \eqref{eq-sigma0}, \eqref{eq-D}, and the bound for \(H\),
\[
 |U(z)|\le
 \frac{C e^{(\pi/2-\eta)|z|^2}}{1+|z|}
 +\frac{C\sqrt{\log(2+|z|)}}{1+|z|}e^{\beta'|z|^2}
 \le C' e^{(\pi/2-\eta)|z|^2}
\]
off the disks. For a point \(z\) inside a disk, every point on its
boundary has modulus at most \(|z|+2d_H\). Applying the maximum
modulus principle gives
\(|U(z)|\le C' e^{(\pi/2-\eta)(|z|+2d_H)^2}\).
After taking logarithms, dividing by \(|z|^2\), and passing to
the upper limit, we get
\(\tau_2(U)\le\pi/2-\eta=\pi\varepsilon/2\). Letting
\(\varepsilon\downarrow2\beta/\pi\) gives \(\tau_2(U)\le\tau_2(f)\).
\end{proof}

\subsection{Translation of the orthogonality condition}

For \(a\in\C\), let \(W_a\) be the Weyl operator
\begin{equation}\label{eq-translation}
 (W_a h)(z)=k_a(z)h(z-a),\qquad h\in\F_\pi^2.
\end{equation}
The defining formula for \(k_a\) gives
\[
 |k_a(z)|^2 e^{-\pi|z|^2}
 =e^{2\pi\operatorname{Re}(z\overline a)-\pi|a|^2-\pi|z|^2}
 =e^{-\pi|z-a|^2}.
\]
Changing variables in the norm integral gives
\(\|W_a h\|_\pi=\|h\|_\pi\). The identity
\(k_{-a}(z)k_a(z+a)=1\) gives \(W_{-a}W_a=I\), so \(W_a\) is
unitary.
Put \(f_a=W_a f\), \(g_a=W_a g\), and \(V_a(z)=V(z-a)\). Since
\(W_a(pf)(z)=p(z-a)f_a(z)\) and translation maps the polynomials
onto themselves,
\begin{equation}\label{eq-covariance1}
 g_a=f_aV_a,\qquad [f_a]_\pi=W_a[f]_\pi,\qquad g_a\perp[f_a]_\pi.
\end{equation}
For a nonzero entire function \(h\) of finite quadratic type and
\(\beta>\tau_2(h)\), we have
\(|h(z-a)|\le C_\beta e^{\beta(|z|+|a|)^2}\).
Multiplication by \(k_a\) adds at most \(O_a(|z|+1)\) to the
logarithm of this bound. Therefore both operations give type at
most \(\beta\), and then at most \(\tau_2(h)\) by decreasing
\(\beta\). Applying the inverse operations gives equality.
In particular \(\tau_2(f_a)=\tau_2(f)\) and
\(\tau_2(V_a)=\tau_2(V)\).
Let \(U_a,S_a\) be the functions in \eqref{eq-US} formed from
\((f_a,g_a,V_a)\). By Proposition~\ref{prop-entire},
\begin{equation}\label{eq-translatedtypes}
 U_a\not\equiv0,\qquad \tau_2(U_a)\le\tau_2(f).
\end{equation}

\begin{lemma}\label{lem-translations}
For every \(a\in\C\),
\begin{equation}\label{eq-covariance2}
 C_{f_a,g_a}(z)=C_{f,g}(z-a),\qquad
 A_{g_a}(z)=A_g(z-a),\qquad S_a(z)=S(z-a).
\end{equation}
\end{lemma}
\begin{proof}
Using \(|k_a(w)|^2e^{-\pi|w|^2}=e^{-\pi|w-a|^2}\) and substituting
\(w=\zeta+a\) gives
\begin{align*}
 C_{f_a,g_a}(z)
 &=\int_\C\frac{f(w-a)\overline{g(w-a)}}{z-w}
           |k_a(w)|^2e^{-\pi|w|^2}\dA(w)\\
 &=\int_\C\frac{f(\zeta)\overline{g(\zeta)}}{z-a-\zeta}
          e^{-\pi|\zeta|^2}\dA(\zeta)=C_{f,g}(z-a).
\end{align*}
The formula for \(A_{g_a}\) follows by taking \(f=g\). The definition
of \(S_a\) then gives
\[
 S_a(z)=V(z-a)C_{f,g}(z-a)-A_g(z-a)=S(z-a).
\]
\end{proof}

\subsection{Zero counting and the contradiction}

\begin{proof}[Proof of Theorem~\ref{thm-main}]
Under \eqref{eq-defect}, Proposition~\ref{prop-polynomial} gives
\(S\in\C[z]\). By Lemma~\ref{lem-translations}, every \(S_a\) is
then a polynomial.

Fix a positive integer \(l\), set \(a_j=j/l\) for \(0\le j<l\), and
choose \(0<\rho<1/(3l)\). By Proposition~\ref{prop-entire},
\[
 I_{g_{a_j}}=\frac{U_{a_j}V_{a_j}}{\sigma_0}-S_{a_j}.
\]
Lemma~\ref{lem-BB}(ii) supplies a lattice thin set
\(\Gamma_j\subset\Lambda_0\) such that \(U_{a_j}V_{a_j}\) has a zero
in each \(D(\lambda,\rho)\) with
\(\lambda\in\Lambda_0\setminus\Gamma_j\). Define
\begin{equation}\label{eq-H}
 G(z)=V(z)\prod_{j=0}^{l-1}U_{a_j}(z+a_j).
\end{equation}
This is a nonzero entire function by \eqref{eq-translatedtypes} and
\(V\not\equiv0\). Since \(V_{a_j}(z+a_j)=V(z)\) for each \(j\), we can
write \(G\) in the form
\[
 G(z)=\bigl(U_{a_j}V_{a_j}\bigr)(z+a_j)
           \prod_{\substack{0\le k<l\\ k\ne j}}U_{a_k}(z+a_k).
\]
It follows that \(G\) has a zero in each disk \(D(\lambda-a_j,\rho)\)
with \(\lambda\in\Lambda_0\setminus\Gamma_j\). The multiplier \(V\)
occurs only once in the product \eqref{eq-H}.

Put
\[
 \Gamma^{(l)}=\bigcup_{j=0}^{l-1}(\Lambda-a_j)=\tfrac1l\Z+i\Z,\qquad
 \mathcal E=\bigcup_{j=0}^{l-1}\bigl((\Gamma_j\cup\{0\})-a_j\bigr).
\]
Every \(w\in\Gamma^{(l)}\setminus\mathcal E\) is the center of a disk
\(D(w,\rho)\) containing a zero of \(G\). Distinct centers are at
distance at least \(1/l>2\rho\), so these zeros are distinct.
If \(|w|<r-\rho\), the whole disk \(D(w,\rho)\) lies in \(D(0,r)\).

The lattice \(\Gamma^{(l)}\) has rectangular cells of area \(1/l\),
centered at the lattice points. The union of the cells centered
in \(D(0,t)\) contains \(D(0,t-1)\) and lies in \(D(0,t+1)\) for
\(t>1\), because every point of a cell is within distance \(1\) of its
center. Comparing areas gives
\[
 \#(\Gamma^{(l)}\cap D(0,t))=l\pi t^2+O_l(t+1).
\]
Thus for \(r>\rho\),
\begin{equation}\label{eq-zerobound}
 n_G(r)\ge\#\bigl(\Gamma^{(l)}\cap D(0,r-\rho)\bigr)-N_{\mathcal E}(r)
 \ge l\pi r^2-N_{\mathcal E}(r)-O_{l,\rho}(r+1)
\end{equation}

Choose \(r_k\to\infty\) as in Lemma~\ref{lem-commonradii}.
For fixed \(0<\delta<1\), integrate \eqref{eq-zerobound} and use
\(N_{\mathcal E}(r)\le N_{\mathcal E}(r_k)\) for \(r\le r_k\) to obtain
\begin{align}
 \int_{\delta r_k}^{r_k}\frac{n_G(r)}{r}\,dr
 &\ge\frac{l\pi}{2}(1-\delta^2)r_k^2
      -N_{\mathcal E}(r_k)\log\frac1\delta
      -O_{l,\rho,\delta}(r_k+1)\notag\\
 &=\frac{l\pi}{2}(1-\delta^2)r_k^2-o(r_k^2).
 \label{eq-jensenlower}
\end{align}
The exceptional count only enters at the outer radius \(r_k\).
This is why one common sequence of radii from
Lemma~\ref{lem-commonradii} suffices.

Write \(G(z)=z^qG_0(z)\), where \(G_0(0)\ne0\). Jensen's formula
\cite[Chapter~5, Section~3]{Ahlfors} expresses the integral as a
difference of circular means of \(\log|G|\).
The circular mean at \(\delta r_k\) is at least
\(q\log(\delta r_k)+\log|G_0(0)|\) by the submean inequality for
\(\log|G_0|\). Hence
\begin{align*}
 \int_{\delta r_k}^{r_k}\frac{n_G(r)}r\,dr
 &=\frac1{2\pi}\int_0^{2\pi}
  \bigl(\log|G(r_ke^{i\theta})|
        -\log|G(\delta r_ke^{i\theta})|\bigr)\,d\theta\\
 &\le\log M_G(r_k)-q\log(\delta r_k)-\log|G_0(0)|.
\end{align*}
Dividing by \(r_k^2\) and taking the upper limit in conjunction with
\eqref{eq-jensenlower} gives
\[
 \tfrac{l\pi}{2}(1-\delta^2)
 \le\limsup_{k\to\infty}\frac{\log M_G(r_k)}{r_k^2}
 \le\tau_2(G).
\]
Letting \(\delta\downarrow0\) yields \(l\pi/2\le\tau_2(G)\).

The definition \eqref{eq-H} gives
\[
 M_G(r)\le M_V(r)\prod_{j=0}^{l-1}M_{U_{a_j}}(r+|a_j|).
\]
Taking logarithms, dividing by \(r^2\), and using
\((r+|a_j|)^2/r^2\to1\), we obtain from \eqref{eq-translatedtypes}
\[
 \tau_2(G)\le\tau_2(V)+\sum_{j=0}^{l-1}\tau_2(U_{a_j})
 \le\tau_2(V)+l\tau_2(f).
\]
Combining the bounds gives
\[
 l\bigl(\pi/2-\tau_2(f)\bigr)\le\tau_2(V).
\]
The two quantities \(\pi/2-\tau_2(f)>0\) and \(\tau_2(V)<\infty\)
are fixed by the original pair \(f,V\). Taking \(l\) arbitrarily
large gives a contradiction.
\end{proof}

\section{Quadratic comparison principle and small-value theorem}\label{sec-comparison}

\subsection{A quadratic comparison principle}

\begin{lemma}\label{lem-comparison}
Let \(u\ge0\) and \(v\not\equiv-\infty\) be subharmonic on \(\C\). Suppose
that for some \(C>0\), \(b\in\R\) and \(0\le\beta<\pi/2\),
\[
 u(z)\le C(1+|z|^2),\qquad v(z)\le\beta|z|^2+b,\qquad
 \Delta v\ge 2\pi\mathbf 1_{\{u>0\}}\,dA.
\]
The last inequality is understood in the sense of measures.
Then \(u\equiv0\).
\end{lemma}

\begin{proof}
Choose \(\eta,\kappa>0\) such that
\[
 \eta C<\kappa<\frac{\pi}{2}-\beta.
\]
The proof proceeds by adding a small multiple of \(u\) at each step.
We show that
\begin{equation}\label{eq-induction}
 v+k\eta u\le\beta|z|^2+b,\qquad k=0,1,2,\ldots.
\end{equation}
The case \(k=0\) is the assumed upper bound for \(v\).
Assume that \eqref{eq-induction} holds for \(k\), and put
\[
 f=v+(k+1)\eta u,\qquad h(z)=f(z)-\beta|z|^2-b.
\]
Then \(h\le\eta u\le\eta C(1+|z|^2)\).
Suppose that \(h(z_1)>0\) for some \(z_1\in\C\). Since \(\eta C<\kappa\),
the upper semicontinuous function \(h(z)-\kappa|z-z_1|^2\)
tends to \(-\infty\) at infinity. It therefore attains
its maximum at a point \(z_0\), and the maximum satisfies
\[
 M=h(z_0)-\kappa|z_0-z_1|^2\ge h(z_1)>0.
\]
Equivalently, the quadratic function
\[
 \phi(z)=\beta|z|^2+b+\kappa|z-z_1|^2+M
\]
satisfies \(f\le\phi\) on \(\C\) and \(f(z_0)=\phi(z_0)\).
On \(Z=\{u=0\}\) we have \(f=v\le\beta|z|^2+b\), so
\(\phi-f\ge M\) there. By the submean inequality,
\(\int_{D(z_0,r)}(f-f(z_0))\dA\ge0\). Using the equality at \(z_0\),
we get
\begin{align*}
 M|Z\cap D(z_0,r)|
 &\le \int_{D(z_0,r)}(\phi-f)\dA\\
 &= \int_{D(z_0,r)}(\phi-\phi(z_0))\dA
    -\int_{D(z_0,r)}(f-f(z_0))\dA\\
 &\le \int_{D(z_0,r)}(\phi-\phi(z_0))\dA
 =\frac{\pi(\beta+\kappa)}2r^4.
\end{align*}
The last equality follows by writing \(z=z_0+\zeta\) and integrating
the quadratic term \((\beta+\kappa)|\zeta|^2\). The linear terms
integrate to zero. Since
\(\log(r/|z-z_0|)=\int_{|z-z_0|}^r ds/s\) in \(D(z_0,r)\),
Tonelli's theorem gives
\[
 \int_{Z\cap D(z_0,r)}\log\frac{r}{|z-z_0|}\dA(z)
 = \int_0^r\frac{|Z\cap D(z_0,s)|}{s}\,ds
 \le \frac{\pi(\beta+\kappa)}{8M}r^4.
\]
The subharmonic Jensen formula, together with
\(\Delta f\ge\Delta v\ge 2\pi\mathbf 1_{\C\setminus Z}\,dA\), gives
\begin{align*}
 \frac{1}{2\pi}\int_0^{2\pi} f(z_0+re^{i\theta})\,d\theta-f(z_0)
 &=\frac{1}{2\pi}
   \int_{D(z_0,r)}\log\frac r{|z-z_0|}\,d(\Delta f)(z)\\
 &\ge \int_{D(z_0,r)\setminus Z}\log\frac r{|z-z_0|}\dA(z)\\
 &\ge \frac\pi2r^2-\frac{\pi(\beta+\kappa)}{8M}r^4.
\end{align*}
In the last line we use the integral of the logarithm over the
whole disk,
\[
 \int_{D(z_0,r)}\log\frac r{|z-z_0|}\dA(z)
   =2\pi\int_0^r s\log\frac rs\,ds=\frac\pi2r^2.
\]
On the other hand, since \(f\le\phi\) and \(f(z_0)=\phi(z_0)\), the
same circular mean difference is at most
\[
 \frac{1}{2\pi}\int_0^{2\pi}
 \bigl(\phi(z_0+re^{i\theta})-\phi(z_0)\bigr)\,d\theta
 =(\beta+\kappa)r^2.
\]
Comparing these two bounds and letting \(r\downarrow0\) gives
\(\beta+\kappa\ge\pi/2\), a contradiction. We conclude that \(h\le0\),
which proves \eqref{eq-induction} for \(k+1\).

Letting \(k\to\infty\) in \eqref{eq-induction}, we obtain \(u=0\) wherever
\(v\) is finite. Since \(v\) is finite almost everywhere, the submean
inequality and \(u\ge0\) give \(u\equiv0\).
\end{proof}

\begin{remark}\label{rem-strictness}
The restriction \(\beta<\pi/2\) is sharp. The functions \(u=1\) and
\(v(z)=(\pi/2)|z|^2\) satisfy \(\Delta v=2\pi\,dA\).
\end{remark}

\subsection{A small-value theorem}\label{s4-small-value-section}

For \(w\in\Lambda\), write
\[
 Q_w=w+[-1/2,1/2)+i[-1/2,1/2),\qquad
 r_w=\bigl(\log(2+|w|)\bigr)^{-2}.
\]
In particular, \(r_w<1/4\) for all sufficiently large \(|w|\).

\begin{theorem}\label{th-small}
Let \(U\) be a nonzero entire function with \(\tau_2(U)<\pi/2\), and
let \(S\) be an entire function with \(\tau_2(S)<\infty\). Set
\[
 E_U=\{w\in\Lambda:D(w,r_w)\cap Z(U)\ne\varnothing\}.
\]
Suppose \(E\subset\Lambda\) has zero radial density and there are constants
\(C>0\) and \(N\ge0\) such that, for every sufficiently large
\(w\in\Lambda\setminus(E_U\cup E)\),
\begin{equation}\label{eq:smallvalues}
 \min_{|z-w|=s}|S(z)|\le C(1+|w|)^N,\qquad r_w/2\le s\le r_w.
\end{equation}
Then \(S\) is a polynomial.
\end{theorem}

\begin{lemma}\label{lem:local}
There is an absolute constant \(C_0>0\) with the following property.
Let \(w\in\Lambda\), \(0<t<1/4\), and let \(q\) be a nonnegative
continuous function in \(W^{1,2}(D(w,2))\). If
\(\min_{|z-w|=s}q(z)\le\epsilon\) for every \(t/2\le s\le t\), where
\(\epsilon\ge0\), then
\begin{equation}\label{eq:local}
 \int_{Q_w}q^2\dA
 \le C_0\left(1+\log\frac1t\right)
         \int_{D(w,1)}|\nabla q|^2\dA+C_0\epsilon^2.
\end{equation}
\end{lemma}

\begin{proof}
Set \(q_*=(q-\epsilon)_+\). The Sobolev chain rule gives
\(|\nabla q_*|\le|\nabla q|\) almost everywhere. Since \(|Q_w|=1\),
\[
 \int_{Q_w}q^2\dA\le2\int_{Q_w}q_*^2\dA+2\epsilon^2.
\]
Also \(q_*\) has a zero on each circle \(|z-w|=s\) with
\(t/2\le s\le t\). It suffices to prove the estimate for \(\epsilon=0\).
After translating to \(w=0\), set
\(\mathcal D=\int_{D(0,1)}|\nabla q|^2\dA\) and write
\[
 q_r(\theta)=q(re^{i\theta}),\qquad
 m(r)=\frac{1}{2\pi}\int_0^{2\pi}q_r(\theta)\,d\theta,\qquad
 A(r)=\int_0^{2\pi}|q_r'(\theta)|^2\,d\theta.
\]
On compact annuli within \(D(0,1)\), Sobolev slicing in polar
coordinates gives \(q_r\in W^{1,2}\) for almost every \(r\), along
with local absolute continuity of \(m\). For almost every \(r\),
\[
 m'(r)=\frac{1}{2\pi}\int_0^{2\pi}
       \partial_rq(re^{i\theta})\,d\theta,\qquad
 |m'(r)|^2\le\frac{1}{2\pi}\int_0^{2\pi}
                       |\partial_rq(re^{i\theta})|^2\,d\theta.
\]
The polar-coordinate identity for the energy is
\[
 \mathcal D
   =\int_0^1\int_0^{2\pi}|\partial_rq(re^{i\theta})|^2\,d\theta\,r\,dr
       +\int_0^1 A(r)\,\frac{dr}{r}.
\]
In particular,
\(\int_0^1 A(r)\,dr/r\le\mathcal D\) and
\(2\pi\int_0^1|m'(r)|^2r\,dr\le\mathcal D\).

The interval \([t/2,t]\) has measure \(\log2\) with respect to \(dr/r\).
We may choose \(s_0\) in this interval with
\(A(s_0)\le\mathcal D/\log2\) for which \(q_{s_0}\in W^{1,2}\).
Its continuous representative vanishes at some \(\theta_0\) by the
hypothesis. Hence
\[
 |q_{s_0}(\theta)|^2
 =\left|\int_{\theta_0}^{\theta}q_{s_0}'(x)\,dx\right|^2
 \le 2\pi A(s_0)
\]
for all \(\theta\in[0,2\pi]\), so
\(|m(s_0)|^2\le (2\pi/\log2)\mathcal D\).
For any \(0<r<1\), let \(J\) be the interval with endpoints \(r\) and \(s_0\).
Using Cauchy--Schwarz with weights \(x\) and \(1/x\), we obtain
\begin{align*}
 |m(r)-m(s_0)|^2
 &\le \left(\int_J |m'(x)|^2x\,dx\right)
        \left(\int_J \frac{dx}{x}\right)\\
 &\le \frac{\mathcal D}{2\pi}\left|\log\frac r{s_0}\right|.
\end{align*}
Since
\[
 \int_0^1 r\left|\log\frac r{s_0}\right|\,dr
 =\frac12\log\frac1{s_0}-\frac14+\frac{s_0^2}{2}
 \le\frac12\log\frac2t,
\]
it follows that
\[
 \int_0^1 |m(r)|^2 r\,dr
 \le |m(s_0)|^2+\frac{\mathcal D}{2\pi}\log\frac2t
 \le C\left(1+\log\frac1t\right)\mathcal D.
\]
The Poincar\'e inequality on the circle gives
\(\int_0^{2\pi}|q_r-m(r)|^2\,d\theta\le A(r)\). Because
\(q_r-m(r)\) has mean zero and \(Q_0\subset D(0,1)\),
\begin{align*}
 \int_{Q_0}q^2\dA
 &\le\int_0^1\int_0^{2\pi}
       \bigl(|q_r(\theta)-m(r)|^2+|m(r)|^2\bigr)\,d\theta\,r\,dr\\
 &\le \int_0^1 A(r)\frac{dr}{r} +2\pi\int_0^1|m(r)|^2r\,dr\\
 &\le C\left(1+\log\frac1t\right)\mathcal D,
\end{align*}
where we used \(r\le 1/r\) for \(r\in(0,1)\).
The truncation estimate then gives \eqref{eq:local}.
\end{proof}

\begin{lemma}\label{lem:qDeltaq}
If \(q\) is nonnegative and subharmonic on \(\C\), then
\(q\in W^{1,2}_{\mathrm{loc}}(\C)\). For every real
\(\chi\in C_c^\infty(\C)\),
\begin{equation}\label{eq:caccioppoli}
 \int\chi^2|\nabla q|^2\dA\le 4\int q^2|\nabla\chi|^2\dA.
\end{equation}
\end{lemma}

\begin{proof}
First suppose \(q\) is smooth. Integration by parts gives
\[
 0\le\int\chi^2q\,\Delta q\dA
 =-\int\chi^2|\nabla q|^2\dA
   -2\int\chi q\,\nabla\chi\cdot\nabla q\dA.
\]
Write \(A=\int\chi^2|\nabla q|^2\dA\) and
\(B=\int q^2|\nabla\chi|^2\dA\).
Taking absolute values in the mixed term and applying
Cauchy--Schwarz gives \(A\le2\sqrt{AB}\). If \(A>0\), divide by
\(\sqrt A\) and square to obtain \(A\le4B\). This inequality also
holds when \(A=0\).

For a general nonnegative subharmonic \(q\), choose a smooth,
nonnegative approximate identity \(\eta_\varepsilon\) and set
\(q_\varepsilon=q*\eta_\varepsilon\). Then \(q_\varepsilon\ge0\) and
\(\Delta q_\varepsilon=(\Delta q)*\eta_\varepsilon\ge0\).
Since \(q\) is locally bounded, \(q_\varepsilon\to q\) in
\(L^2_{\mathrm{loc}}\). Taking cutoffs that equal one on compact sets
in the smooth inequality gives a locally uniform \(L^2\) bound on
\(\nabla q_\varepsilon\). On each compact disk the weak limits of
these gradients must be \(\nabla q\) in the distributional sense.
Thus \(q\in W^{1,2}_{\mathrm{loc}}\), and for any fixed \(\chi\),
\[
 \int\chi^2|\nabla q|^2\dA
 \le\liminf_{\varepsilon\downarrow0}
         \int\chi^2|\nabla q_\varepsilon|^2\dA
 \le4\int q^2|\nabla\chi|^2\dA.
\]
Here the first inequality uses weak lower semicontinuity and the
second uses strong local \(L^2\) convergence on the support of \(\chi\).
\end{proof}

\begin{proof}[Proof of Theorem~\ref{th-small}]
Suppose \(S\) is transcendental, and set \(L(R)=\log^+M_S(R)\).
For every nonzero Taylor coefficient \(a_k\), Cauchy's inequality
gives \(L(R)\ge k\log R+\log|a_k|\). Since the indices of these
coefficients are unbounded, \(L(R)/\log R\to\infty\).
Finite quadratic type gives \(L(R)=O(R^2)\). In particular, the
supremum below is finite and positive, so we may choose \(R_j\ge j\)
such that
\[
 \frac{L(R_j)}{R_j^2}\ge\frac12\sup_{R\ge j}\frac{L(R)}{R^2}>0.
\]
Put
\[
 u_j(z)=\frac{\log^+|S(R_jz)|}{L(R_j)}.
\]
The functions \(u_j\) are continuous, nonnegative and subharmonic.
For \(|z|\le1\), monotonicity of \(M_S\) gives \(u_j(z)\le1\),
with equality at some \(z_j\) on the unit circle. For \(|z|>1\),
the choice of \(R_j\) gives
\[
 u_j(z)\le\frac{L(R_j|z|)}{L(R_j)}
 = |z|^2\,\frac{L(R_j|z|)/(R_j|z|)^2}{L(R_j)/R_j^2}
 \le2|z|^2.
\]
Thus \(0\le u_j(z)\le2(1+|z|^2)\).
By compactness for subharmonic functions \cite{Ransford}, a
subsequence converges in \(L^1_{\mathrm{loc}}\) to a subharmonic
function \(u\). The alternative of local uniform convergence to
\(-\infty\) is excluded by \(u_j\ge0\). The bounds pass to the limit
almost everywhere and then everywhere for the subharmonic
representative. Indeed, its value at a point is the limit of its
disk means as the radii decrease to zero. On every compact disk \(K\)
we have
\[
 \int_K|u_j-u|^2\dA\le C_K\int_K|u_j-u|\dA\longrightarrow0.
\]
The lower bound at \(z_j\) persists in the limit by the submean
inequality, since \(D(z_j,1)\subset D(0,2)\). Hence
\[
 \int_{D(0,2)}u\dA
 =\lim_{j\to\infty}\int_{D(0,2)}u_j\dA
 \ge\liminf_{j\to\infty}\int_{D(z_j,1)}u_j\dA\ge\pi.
\]

Choose \(\beta\) with \(\tau_2(U)<\beta<\pi/2\). For
\(v_j(z)=R_j^{-2}\log|U(R_jz)|\), the bound
\(\log|U(\zeta)|\le\beta|\zeta|^2+b\), \(b\in\R\), gives
\[
 v_j(z)\le\beta|z|^2+b/R_j^2.
\]
Fix \(z_*\) with \(U(z_*)\ne0\). The points \(z_*/R_j\) tend to
\(0\), and \(v_j(z_*/R_j)=R_j^{-2}\log|U(z_*)|\to0\).
Consequently, the sequence cannot converge locally uniformly to
\(-\infty\). Subharmonic compactness gives, after another subsequence,
\(v_j\to v\) in \(L^1_{\mathrm{loc}}\), where \(v\not\equiv-\infty\)
is subharmonic. Taking disk averages of the upper bound gives
\(v(z)\le\beta|z|^2\).

Define \(\chi_j=\sum_{w\in E_U}\mathbf 1_{Q_w/R_j}\).
The lattice cells are disjoint, so \(0\le \chi_j\le1\).
Weak-* compactness on disks and a diagonal subsequence give a
function \(0\le\rho\le1\) with
\(\int\chi_j p\,dA\to\int\rho p\,dA\) for every compactly supported
\(p\in L^1(\C)\).
For every sufficiently large \(w\in E_U\), choose a zero \(\lambda_w\)
of \(U\) in \(D(w,r_w)\). These disks are disjoint. The distributional identity for the logarithm of a holomorphic
function gives
\[
 \Delta v_j=\frac{2\pi}{R_j^2}
             \sum_{\lambda\in Z(U)}m_U(\lambda)\delta_{\lambda/R_j}.
\]
For each nonnegative \(\varphi\in C_c^\infty(\C)\), compare the mass
at the selected zero with the area \(R_j^{-2}\) of its cell.
For \(\zeta\in Q_w/R_j\) and a selected \(\lambda_w\),
\(|\zeta-\lambda_w/R_j|\le\sqrt2/R_j\), so
\begin{align*}
 \left|\frac{\varphi(\lambda_w/R_j)}{R_j^2}
            -\int_{Q_w/R_j}\varphi(\zeta)\dA(\zeta)\right|
 &\le\int_{Q_w/R_j}|\varphi(\lambda_w/R_j)-\varphi(\zeta)|\dA(\zeta)\\
 &\le\frac{\sqrt2\|\nabla\varphi\|_\infty}{R_j^3}.
\end{align*}
Only \(O_\varphi(R_j^2)\) cells can meet the support of \(\varphi\)
or have a selected zero there. The finitely many cells with no
selected zero have total area \(O(R_j^{-2})\). Since the selected
zeros are distinct, summing the cell estimates gives
\[
 \int\varphi\,d(\Delta v_j)
 \ge\frac{2\pi}{R_j^2}\sum_{\substack{w\in E_U\\ |w|\ \mathrm{large}}}
      \varphi(\lambda_w/R_j)
 \ge 2\pi\int\varphi\chi_j\dA-O_\varphi(R_j^{-1}).
\]
On the left,
\(\int\varphi\,d(\Delta v_j)=\int v_j\Delta\varphi\dA\to
\int v\Delta\varphi\dA\) by \(L^1_{\mathrm{loc}}\) convergence.
On the right, the weak-* limit gives
\(2\pi\int\varphi\chi_j\dA\to2\pi\int\varphi\rho\dA\).
The resulting inequality for nonnegative test functions is
\begin{equation}\label{eq:laptheta}
 \Delta v\ge 2\pi\rho\,dA.
\end{equation}

Fix \(r>0\). For large \(j\), if \(Q_w/R_j\) meets \(D(0,r)\),
then \(|w|<(r+1)R_j\). Suppose such a \(w\) is sufficiently large
and belongs to \(\Lambda\setminus(E_U\cup E)\). By
\eqref{eq:smallvalues}, the function \(q_j(\zeta)=u_j(\zeta/R_j)\)
has minimum at most
\[
 \epsilon_j=\frac{\log^+C+N\log(1+(r+1)R_j)}{L(R_j)}
\]
on each circle \(|\zeta-w|=s\), \(r_w/2\le s\le r_w\).
Since \(L(R)/\log R\to\infty\), we have \(\epsilon_j\to0\).

We apply Lemma~\ref{lem:local} with \(q=q_j\), \(t=r_w\).
Under \(\zeta=R_jz\), area is multiplied by \(R_j^2\) and
\(\nabla q_j(R_j z)=R_j^{-1}\nabla u_j(z)\). Hence
\[
 \int_{Q_w} q_j^2\dA=R_j^2\int_{Q_w/R_j}u_j^2\dA,\qquad
 \int_{D(w,1)}|\nabla q_j|^2\dA
   =\int_{D(w/R_j,1/R_j)}|\nabla u_j|^2\dA.
\]
Writing \(1+\log(1/r_w)\le C\log(2/r_w)\), the local estimate becomes
\[
 \int_{Q_w/R_j}u_j^2\dA
 \le\frac{C\log(2/r_w)}{R_j^2}
          \int_{D(w/R_j,1/R_j)}|\nabla u_j|^2\dA
       +\frac{C\epsilon_j^2}{R_j^2}.
\]
There are \(O_r(R_j^2)\) cells under consideration. Moreover,
the disks \(D(w/R_j,1/R_j)\) have bounded overlap, since any point
in such a disk lies within distance one of \(w\) after rescaling.
These disks lie in \(D(0,r+2)\) for large \(j\). A cutoff equal
to one on \(D(0,r+2)\), together with Lemma~\ref{lem:qDeltaq}
and \(u_j(z)\le2(1+|z|^2)\), gives
\[
 \int_{D(0,r+2)}|\nabla u_j|^2\dA\le C_r.
\]
We also have \(\log(2/r_w)=O_r(\log\log(3+R_j))\) when
\(|w|<(r+1)R_j\). Summing the estimate on the cells under consideration
therefore gives an upper bound
\[
 C_r\frac{\log\log(3+R_j)}{R_j^2}+C_r\epsilon_j^2
\]
for the contribution of these cells to
\(\int_{D(0,r)}u_j^2(1-\chi_j)\dA\).
The remaining cells contributing to this integral come from \(E\)
or from the finitely many small \(w\) excluded from
\eqref{eq:smallvalues}. Their total area in \(D(0,r)\) is at most
\(R_j^{-2}N_E((r+1)R_j)+O(R_j^{-2})=o(1)\), and \(u_j\) is uniformly
bounded there. We conclude that
\[
 \int_{D(0,r)} u_j^2(1-\chi_j)\dA\longrightarrow0.
\]
Strong \(L^2\) convergence on \(D(0,r)\) gives
\[
 \|u_j^2-u^2\|_{L^1(D(0,r))}
 \le\|u_j-u\|_{L^2(D(0,r))}\|u_j+u\|_{L^2(D(0,r))}
 \longrightarrow0.
\]
We may replace \(u_j^2\) by \(u^2\) in the last integral since
\(0\le1-\chi_j\le1\). Weak-* convergence of \(\chi_j\), tested
against \(u^2\mathbf 1_{D(0,r)}\), now gives
\(\int_{D(0,r)}u^2(1-\rho)\dA=0\).
As \(r\) was arbitrary, \(\rho=1\) almost everywhere on \(\{u>0\}\).
Thus \eqref{eq:laptheta} gives
\(\Delta v\ge 2\pi\mathbf 1_{\{u>0\}}\,dA\). Lemma~\ref{lem-comparison}
implies \(u\equiv0\), contradicting the lower bound for its integral
on \(D(0,2)\).
\end{proof}

\subsection{Polynomiality}\label{s6-polynomial-section}

\begin{lemma}\label{lem-jensen}
Let \(H\) be a nonzero entire function of finite quadratic type.
Then \(n_H(r)=O(r^2)\), and for each fixed \(d\ge0\),
\begin{equation}\label{eq-jensen-bound}
 \limsup_{R\to\infty}\frac1{\log R}
 \int_1^R\frac{n_H(t+d)}{t^3}\,dt\le2\tau_2(H).
\end{equation}
\end{lemma}

\begin{proof}
Let \(J_H(R)=\int_1^R n_H(t)\,dt/t\).
Jensen's formula \cite{Ransford} gives
\[
 J_H(R)=\frac{1}{2\pi}\int_0^{2\pi}
       \bigl(\log|H(Re^{i\theta})|-\log|H(e^{i\theta})|\bigr)\,d\theta.
\]
This identity also counts a zero at the origin, since a factor
\(z^m\) contributes \(m\log R\) on both sides.
By the definition of quadratic type,
\[
 J_H(R)\le \log M_H(R)+O_H(1)\le\tau_2(H)R^2+o(R^2).
\]
Monotonicity of \(n_H\) gives, for \(r\ge1\),
\[
 n_H(r)\log2\le\int_r^{2r}\frac{n_H(t)}{t}\,dt
 \le J_H(2r)=O(r^2).
\]
Finally, integration by parts gives
\[
 \int_1^R\frac{n_H(t)}{t^3}\,dt
 = \frac{J_H(R)}{R^2}+2\int_1^R\frac{J_H(t)}{t^3}\,dt.
\]
The first term on the right is bounded. Given \(\varepsilon>0\),
we have \(J_H(t)\le(\tau_2(H)+\varepsilon)t^2\) for all sufficiently
large \(t\). Hence the second term is at most
\(2(\tau_2(H)+\varepsilon)\log R+O_{\varepsilon,H}(1)\).
Dividing by \(\log R\) and letting \(\varepsilon\downarrow0\)
proves \eqref{eq-jensen-bound} for \(d=0\).

For fixed \(d>0\), the substitution \(s=t+d\) gives the integral
of \(n_H(s)(s-d)^{-3}\) over \([1+d,R+d]\). Replacing
\((s-d)^{-3}\) by \(s^{-3}\) changes this integral by \(O_d(1)\)
because their difference is \(O_d(s^{-4})\) and \(n_H(s)=O(s^2)\).
Changing the endpoints to \(1\) and \(R\) also changes the integral
by a bounded amount. This proves \eqref{eq-jensen-bound} for \(d>0\).
\end{proof}

\begin{lemma}\label{lem-generic}
Let \(V\) be a nonzero entire function of finite quadratic type. For
almost every \(a\in\C\), the set
\[
 E_a=\{w\in\Lambda:D(w,r_w)\cap Z(V(\cdot-a))\ne\varnothing\}
\]
has zero radial density.
\end{lemma}

\begin{proof}
Fix \(r\ge1\) and let \(R\) be large. We estimate
\(\int_{D(0,r)}N_{E_a}(R)\,dA(a)\) by splitting the count at
\(|w|=\sqrt R\). The inner part contributes at most \(C_rR\),
since there are \(O(R)\) lattice points in \(D(0,\sqrt R)\).
For the other cells, use \(Z(V(\cdot-a))=Z(V)+a\) to write
\[
 w\in E_a
 \quad\Longleftrightarrow\quad
 a\in\bigcup_{\lambda\in Z(V)}D(w-\lambda,r_w).
\]
Thus the contribution from the annulus is at most
\[
 \sum_{\substack{w\in\Lambda\\\sqrt R\le|w|<R}}
 \sum_{\lambda\in Z(V)}
       |D(0,r)\cap D(w-\lambda,r_w)|.
\]
If an intersection is nonempty, then
\(|w-\lambda|<r+r_w<r+1\) and
\(|\lambda|<|w|+r+1<R+r+1\).
For each zero \(\lambda\), the first inequality restricts \(w\)
to at most \(C_r\) lattice points. On the annulus,
\(\log(2+|w|)\ge\tfrac12\log R\), so each intersection has area
at most \(\pi r_w^2\le C(\log R)^{-4}\). The annular sum is
therefore bounded by
\[
 \frac{C_r n_V(R+r+1)}{(\log R)^4}
 \le C_{r,V}\frac{R^2}{(\log R)^4},
\]
where the last inequality uses Lemma~\ref{lem-jensen}.
Adding the inner part gives
\[
 \int_{D(0,r)}N_{E_a}(R)\,dA(a)
 \le C_{r,V}\left(R+\frac{R^2}{(\log R)^4}\right).
\]
At the dyadic radii \(R=2^m\), Tonelli's theorem gives
\begin{align*}
 \int_{D(0,r)}\sum_{m=1}^\infty4^{-m}N_{E_a}(2^m)\,dA(a)
 &=\sum_{m=1}^\infty4^{-m}\int_{D(0,r)}N_{E_a}(2^m)\,dA(a)\\
 &\le C_{r,V}\sum_{m=1}^\infty \bigl(2^{-m}+m^{-4}\bigr)<\infty.
\end{align*}
It follows that the series inside the integral converges for
almost every \(a\in D(0,r)\), so \(4^{-m}N_{E_a}(2^m)\to0\).
For \(2^m\le R<2^{m+1}\), monotonicity gives
\[
 \frac{N_{E_a}(R)}{R^2}
 \le4\,\frac{N_{E_a}(2^{m+1})}{4^{m+1}}\longrightarrow0.
\]
Taking \(r=1,2,\ldots\) proves the assertion for almost every \(a\in\C\).
\end{proof}

\begin{proposition}\label{prop-polynomial}
Under \eqref{eq-defect}, the function \(S\) defined in \eqref{eq-US}
and all its translates \(S_a\) are polynomials.
\end{proposition}

\begin{proof}
Choose \(a\) by Lemma~\ref{lem-generic} so that
\[
 E_{V_a}=\{w\in\Lambda:D(w,r_w)\cap Z(V_a)\ne\varnothing\}
\]
has zero radial density. Let
\(E_{U_a}=\{w\in\Lambda:D(w,r_w)\cap Z(U_a)\ne\varnothing\}\).
For \(w\in\Lambda\setminus(E_{U_a}\cup E_{V_a})\), the product
\(U_aV_a\) is zero-free in \(D(w,r_w)\).

Fix such a \(w\) of sufficiently large modulus and choose
\(r_w/2\le s\le r_w\). On the circle \(|z-w|=s\), the evaluation
bound \eqref{eq-eval}, unitarity of the Weyl operator, and
\eqref{eq-circle} give
\[
 |g_a(z)|e^{-\pi|z|^2/2}\le\|g\|_\pi,\qquad
 |\sigma_0(z)|e^{-\pi|z|^2/2}\ge c\,\frac{s}{1+|w|}.
\]
The bound \eqref{eq-D} gives
\[
 |B_{g_a}(z)|\le C_a\,\frac{\sqrt{\log(2+|w|)}}{1+|w|},
\]
since \(\bigl||z|-|w|\bigr|\le s<1/4\).
By boundedness of \(A_{g_a}\) and \eqref{eq-I},
\begin{align*}
 |I_{g_a}(z)|
 &\le |A_{g_a}(z)|+
      \frac{|g_a(z)|}{|\sigma_0(z)|}|B_{g_a}(z)|\\
 &\le C_a\left(1+\frac{\sqrt{\log(2+|w|)}}{s}\right)\\
 &\le C_a\bigl(1+2(\log(2+|w|))^{5/2}\bigr)
 \le C_a'(1+|w|)^3.
\end{align*}
The third line uses \(s\ge r_w/2\).

Suppose \(\min_{|z-w|=s}|S_a(z)|>C_a'(1+|w|)^3\).
Both \(U_aV_a\) and \(\sigma_0S_a\) are entire, and \(\sigma_0\)
has no zeros on \(\partial D(w,s)\).
The translated identity \eqref{eq-identity} gives
\[
 |U_aV_a-\sigma_0S_a|=|\sigma_0I_{g_a}|<|\sigma_0S_a|
 \qquad \text{on }\partial D(w,s).
\]
By Rouch\'e's theorem \cite{Ahlfors}, the two entire functions
have the same number of zeros in \(D(w,s)\), with multiplicity.
The product \(\sigma_0S_a\) vanishes at \(w\). This forces a zero
of \(U_aV_a\) in \(D(w,s)\subset D(w,r_w)\), contrary to the
choice of \(w\). Since \(s\) was arbitrary in \([r_w/2,r_w]\),
we have
\[
 \min_{|z-w|=s}|S_a(z)|\le C_a'(1+|w|)^3
 \qquad(r_w/2\le s\le r_w).
\]
Thus \(S_a\) satisfies \eqref{eq:smallvalues} with
\(U=U_a\), \(E=E_{V_a}\), and \(N=3\).

By Proposition~\ref{prop-entire}, \(U_a\not\equiv0\),
\(\tau_2(U_a)<\pi/2\), and \(S_a\) has finite quadratic growth.
Theorem~\ref{th-small} gives \(S_a\in\C[z]\).
Finally, \eqref{eq-covariance2} gives \(S(z)=S_a(z+a)\), and
then every translate is a polynomial.
\end{proof}

\section{The critical type}\label{sec-critical}

Hou and Ma \cite[Example~6.10]{HouMa} use a related product over
regular \(2n\)-gons with radii \(n+(\log n)/n\) to study exact zero
sets at critical density. Here the displacement \(\sqrt n\) gives a
deficit of order \(r^{3/2}\) in the logarithmic growth, which ensures
exponential integrability. We also compare the multiplier norm with
a radial weighted norm to prove Taylor approximation.

Set
\begin{equation}\label{eq-Rn}
 R_n=n+\sqrt n,\qquad
 \Lambda_*=\{R_ne^{\pi i k/n}: n\ge1,\ 0\le k<2n\},
\end{equation}
and define
\begin{equation}\label{eq-Phi}
 \Phi(r)=\sum_{n\ge1}2n\log^+\frac r{R_n},\qquad
 d_*(z)=\min\{1,\dist(z,\Lambda_*)\}.
\end{equation}
We take \(\Phi(0)=0\).

\begin{proposition}\label{prop-product}
The product defining \(F_*\) in Theorem~\ref{thm-endpoint} converges locally
uniformly to an entire function with \(F_*(0)=1\). Its zeros are simple and
form the set \(\Lambda_*\), with distinct zeros at distance at least \(1\).
Moreover,
\begin{equation}\label{eq-F*asymp}
 |F_*(z)|\asymp e^{\Phi(|z|)}d_*(z),\qquad z\in\C,
\end{equation}
and, as \(r\to\infty\),
\begin{equation}\label{eq-F*type}
 \log M_{F_*}(r)=\Phi(r)+O(1)
 =\frac{r^2}{2}-\frac43r^{3/2}+O(r).
\end{equation}
\end{proposition}

\begin{proof}
Fix \(K>0\) and choose \(n_K\) so that \(R_n\ge2K\) for \(n\ge n_K\).
On the closed disk \(\{|z|\le K\}\) we then have
\[
 \left|\frac z{R_n}\right|^{2n}\le4^{-n},\qquad |z|\le K.
\]
Use the power series of \(\log(1-u)\) for \(|u|<1\). For \(n\ge n_K\)
the resulting analytic branch of \(\log(1-(z/R_n)^{2n})\) satisfies
\[
 \left|\log\!\left(1-(z/R_n)^{2n}\right)\right|
 \le \frac{4^{-n}}{1-4^{-n}}\le 2\,4^{-n}.
\]
The series of tail logarithms therefore converges uniformly on the disk,
and its exponential is the nonvanishing limit of the tail products.
Thus the zeros on this disk are exactly those of the finite initial
product. Since \(K\) is arbitrary, \(F_*\) is entire and \(F_*(0)=1\).
Each factor has simple zeros on its own circle, and the circles have
distinct radii, so the zeros of \(F_*\) are simple. Zeros on different
circles are separated by more than \(1\) because \(R_{n+1}-R_n>1\).
On the same circle their least separation is
\[
 2R_n\sin\frac{\pi}{2n}\ge \frac{2R_n}{n}\ge2.
\]

We remove the radial growth of each factor and bound the resulting
product.
For \(z=re^{i\theta}\) with \(r>0\), put
\(q_n(r)=\exp(-2n|\log(r/R_n)|)\). Then
\begin{equation}\label{eq-factors}
 \frac{|1-(z/R_n)^{2n}|}{\exp(2n\log^+(r/R_n))}
 =|1-q_n(r)e^{\varepsilon_n2ni\theta}|,\qquad \varepsilon_n\in\{-1,1\}.
\end{equation}
For \(r\le R_n\) this identity is immediate with \(\varepsilon_n=1\).
For \(r>R_n\), factor \((r/R_n)^{2n}\) out of the modulus to get
\(\varepsilon_n=-1\). Multiplying the identities gives
\[
 |F_*(re^{i\theta})|e^{-\Phi(r)}
 =\prod_{n\ge1}|1-q_n(r)e^{\varepsilon_n2ni\theta}|.
\]
We have
\begin{equation}\label{eq-Cq}
 C_q:=\sup_{r>0}\sum_{n\ge1}q_n(r)<\infty.
\end{equation}
Indeed, \(q_n(r)\le4^{-n}\) when \(R_n\notin[r/2,2r]\). When
\(r/2\le R_n\le2r\), the inequalities \(R_n\le2n\) and
\(\max\{r,R_n\}\le2R_n\) give
\[
 2n|\log(r/R_n)|
 \ge\frac{2n|r-R_n|}{\max\{r,R_n\}}
 \ge\frac12|r-R_n|.
\]
Since the \(R_n\) are more than one unit apart, each set
\(\{n:k\le|r-R_n|<k+1\}\) has at most four elements. Thus the sum in
\eqref{eq-Cq} satisfies
\[
 \sum_{n\ge1}q_n(r)
 \le \sum_{n\ge1}4^{-n}+4\sum_{k\ge0}e^{-k/2},
 \qquad r>0.
\]

Fix \(\delta=1/8\). We consider separately the radii that stay at least
\(\delta\) from every circle and the radii within \(\delta\) of one
circle. In the first case,
\(q_n(r)\le q_*:=\max\{1/4,e^{-\delta/2}\}<1\) for all \(n\). For
\(0\le q\le q_*\), we have
\[
 -\frac{q}{1-q_*}
 \le \log(1-q)\le \log|1-qe^{i\theta}|
 \le \log(1+q)\le q.
\]
The lower bound follows by integrating \(1/(1-t)\) for \(0\le t\le q\).
Summing over \(n\) and using \eqref{eq-Cq} gives
\[
 e^{-C_q/(1-q_*)}
 \le |F_*(z)|e^{-\Phi(r)}\le e^{C_q}.
\]
Every zero of \(F_*\) is at distance at least \(\delta\) from \(z\)
in this case, so \(\delta\le d_*(z)\le1\). This proves
\eqref{eq-F*asymp} for these radii.

If \(|r-R_{n_0}|<\delta\) for some \(n_0\), then \(n_0\) is unique and
\(|r-R_n|>1-\delta\) for each \(n\ne n_0\).
In particular, \(q_n(r)\le q_*\) still holds for \(n\ne n_0\).
Summing the logarithms of these factors gives uniform constants
\(c_1,C_1>0\) such that
\begin{equation}\label{eq-nn0}
 c_1\le\prod_{n\ne n_0}|1-q_n(r)e^{\varepsilon_n2ni\theta}|\le C_1.
\end{equation}
Write \(R=R_{n_0}\) and choose a closest zero
\(\lambda=Re^{i\theta_0}\) on the circle of radius \(R\).
We may write \(z=re^{i(\theta_0+t)}\) with \(|t|\le\pi/(2n_0)\).
Put \(x=2n_0|\log(r/R)|\). The mean-value theorem for \(\log\) yields
\[
 \frac{2n_0}{\max\{r,R\}}|r-R|
 \le x
 \le \frac{2n_0}{\min\{r,R\}}|r-R|.
\]
Since \(n_0\le R\le2n_0\) and \(R/2\le r\le2R\), this gives
\(\frac12|r-R|\le x\le4|r-R|\le4\delta\). Using
\(xe^{-x}\le 1-e^{-x}\le x\), we obtain
\[
 q_{n_0}(r)=e^{-x}\ge e^{-4\delta},\qquad
 1-q_{n_0}(r)\asymp|r-R|.
\]
The distance to \(\lambda\) has the corresponding radial and angular
terms
\begin{equation}\label{eq:nearzero}
 |z-\lambda|^2=(r-R)^2+4rR\sin^2(t/2)
 \asymp(r-R)^2+R^2t^2.
\end{equation}
Since \(e^{2n_0i\theta_0}=1\) and \(|n_0t|\le\pi/2\), the remaining
factor satisfies
\begin{equation}\label{eq-n0factor}
 \begin{aligned}
 |1-q_{n_0}(r)e^{\varepsilon_{n_0}2n_0i\theta}|^2
 &=(1-q_{n_0}(r))^2+4q_{n_0}(r)\sin^2(n_0t)\\
 &\asymp(r-R)^2+R^2t^2\asymp|z-\lambda|^2.
 \end{aligned}
\end{equation}
The bounds \(|r-R|<\delta\) and \(R|t|\le\pi\) show from
\eqref{eq:nearzero} that \(|z-\lambda|\) is bounded by an absolute
constant. Every zero on a different circle is at distance greater than
\(1-\delta\) from \(z\), so after truncating the distance at \(1\),
we obtain
\begin{equation}\label{eq-asmpd*}
 (1-\delta)\min\{1,|z-\lambda|\}\le d_*(z)\le\min\{1,|z-\lambda|\},
 \qquad d_*(z)\asymp |z-\lambda|.
\end{equation}
Together with \eqref{eq-nn0} and \eqref{eq-n0factor}, this proves
\eqref{eq-F*asymp} for \(z\ne0\). At \(z=0\), both sides equal \(1\).

Let \(N(r)=\#\{n:R_n\le r\}\).
Solving \(n+\sqrt n\le r\) for \(\sqrt n\) and then squaring gives
\[
 N(r)=\left\lfloor r+\frac12-\frac12\sqrt{1+4r}\right\rfloor
 =r-\sqrt r+O(1).
\]
On each interval between consecutive \(R_n\), the terms of \(\Phi\)
with \(R_n>r\) vanish and the remaining terms can be differentiated.
Thus, for \(r\notin\{R_n:n\ge1\}\),
\begin{equation}\label{eq-Phi'}
 \Phi'(r)=\frac1r\sum_{n=1}^{N(r)}2n
 =\frac{N(r)(N(r)+1)}r=r-2\sqrt r+O(1).
\end{equation}
The function \(\Phi\) is continuous and piecewise smooth, and
\(\Phi(1)=0\). Integrating \eqref{eq-Phi'} from \(1\) to \(r\) gives
\[
 \Phi(r)=\int_1^r\bigl(t-2\sqrt t+O(1)\bigr)\,dt
 =\frac{r^2}{2}-\frac43r^{3/2}+O(r).
\]
Since \(F_*(0)=1\) and \(F_*\) has \(2n\) zeros on the circle of
radius \(R_n\), Jensen's formula gives
\[
 \Phi(r)=\frac1{2\pi}\int_0^{2\pi}\log|F_*(re^{i\theta})|\,d\theta
 \le\log M_{F_*}(r).
\]
The upper bound \(\log M_{F_*}(r)\le\Phi(r)+O(1)\) follows from
\eqref{eq-F*asymp} and \(d_*\le1\). Together these bounds prove
\eqref{eq-F*type}.
\end{proof}

Define
\begin{equation}\label{eq-weight}
 w(r)=e^{2\Phi(r)-r^2},\qquad \psi(r)=\log w(r).
\end{equation}
By \eqref{eq-F*type}, \(\psi(r)=-\frac83r^{3/2}+O(r)\) as \(r\to\infty\).
In particular, there is \(c>0\) such that \(w(r)\le e^{-cr^{3/2}}\)
for all sufficiently large \(r\). The upper bound in
\eqref{eq-F*asymp} gives
\begin{equation}\label{eq-exponentials}
 \|e^{\lambda z}F_*\|_1^2
 \le C\int_\C e^{2|\lambda||z|}w(|z|)\dA(z)<\infty,\qquad \lambda\in\C.
\end{equation}
Likewise \(pF_*\in\F_1^2\) for each polynomial \(p\).
The same growth estimate gives the pointwise bound
\begin{equation}\label{eq-critical-growth}
 \text{for every }A>0\text{ there is }C_A>0\text{ such that }\quad
 |F_*(z)|\le C_Ae^{|z|^2/2-A|z|},\qquad z\in\C.
\end{equation}

\begin{lemma}\label{lem-norm}
There are constants \(0<c_0\le C_0<\infty\) such that every entire
function \(V\) satisfies
\begin{equation}\label{eq-norm-equivalence}
 c_0\|V\|_w^2\le\|F_*V\|_1^2\le C_0\|V\|_w^2,
 \qquad
 \|V\|_w^2=\int_\C |V(z)|^2w(|z|)\dA(z).
\end{equation}
The integrals may be infinite.
\end{lemma}

\begin{proof}
The upper bound follows from \eqref{eq-F*asymp} because \(d_*\le1\).
For the lower bound, the same pointwise estimate gives direct control
away from the zeros. Near each zero we compare the weighted integral
on a disk with the integral on a surrounding annulus, where the distance
to the zero set is bounded below.

Set \(\delta=1/8\). For \(r\notin\{R_n:n\ge1\}\),
\begin{equation}\label{s7-psi-derivatives}
 \psi'(r)=\frac{2N(r)(N(r)+1)}r-2r,\qquad
 \psi''(r)=-\frac{2N(r)(N(r)+1)}{r^2}-2.
\end{equation}
For \(r\ge1\) we have \(N(r)\le r\) and \(N(r)+1\le2r\), so
\(|\psi''(r)|\le6\) away from the radii \(R_n\).
Since \(N(r)\) changes from \(n-1\) to \(n\) at \(R_n\), the jump in
the first derivative is
\begin{equation}\label{s7-jumps}
 \psi'(R_n+)-\psi'(R_n-)=\frac{4n}{R_n}\le4.
\end{equation}
Fix \(\lambda\in\Lambda_*\) and write \(R=|\lambda|=R_n\).
With \(a=\psi'(R-)=2n(n-1)/R-2R\), we have \(-2R\le a\le0\).
The interval \([R-2\delta,R+2\delta]\) lies in \((1,\infty)\) and
contains no other \(R_m\). Consequently, for \(r\) in this interval
and \(r\ne R\),
\[
 |\psi'(r)-a|
 \le 6|r-R|+4\mathbf 1_{\{r>R\}}.
\]
Integrating on the corresponding side of \(R\) yields
\begin{equation}\label{eq-psi}
 |\psi(r)-\psi(R)-a(r-R)|\le 3|r-R|^2+4|r-R|,
 \qquad |r-R|\le2\delta.
\end{equation}
For \(z=\lambda+\zeta\) with \(|\zeta|<2\delta\), squaring the moduli and
using \(\bigl||z|-R\bigr|\le|\zeta|\) gives
\begin{equation}\label{eq-D(lambda,2delta)}
 \left||z|-R-\operatorname{Re}\frac{\overline\lambda\zeta}{R}\right|
 =\frac{|\zeta|^2-(|z|-R)^2}{2R}
 \le\frac{|\zeta|^2}{2R}.
\end{equation}
Put \(b_\lambda=a\overline\lambda/(2R)\). Because \(|a|\le2R\), the
error in replacing \(a(|z|-R)\) by \(2\operatorname{Re}(b_\lambda\zeta)\)
is at most \(|\zeta|^2\). Adding the error from \eqref{eq-psi} gives
\begin{equation}\label{s7-local-linearization}
 |\psi(|z|)-\psi(R)-2\operatorname{Re}(b_\lambda(z-\lambda))|
 \le C,\qquad z\in D(\lambda,2\delta).
\end{equation}
Consequently, uniformly in \(\lambda\),
\begin{equation}\label{eq-local-weight}
 w(|z|)\asymp w(R)|e^{b_\lambda(z-\lambda)}|^2,\qquad z\in D(\lambda,2\delta).
\end{equation}

For \(h\) holomorphic on \(D(\lambda,2\delta)\), the circular means of
\(|h|^2\) centered at \(\lambda\) are increasing by subharmonicity.
With \(m_h(s)=(2\pi)^{-1}\int_0^{2\pi}|h(\lambda+se^{i\theta})|^2
\,d\theta\), polar integration gives
\begin{align*}
 \int_{D(\lambda,\delta)}|h|^2\dA
 &=2\pi\int_0^\delta m_h(s)s\,ds\le \pi\delta^2 m_h(\delta),\\
 \int_{D(\lambda,2\delta)\setminus D(\lambda,\delta)}|h|^2\dA
 &=2\pi\int_\delta^{2\delta} m_h(s)s\,ds\ge 3\pi\delta^2 m_h(\delta).
\end{align*}
In particular,
\begin{equation}\label{eq-disk-annulus}
 \int_{D(\lambda,\delta)}|h|^2\dA
 \le\frac13\int_{D(\lambda,2\delta)\setminus D(\lambda,\delta)}|h|^2\dA.
\end{equation}
Apply this inequality to \(h(z)=V(z)e^{b_\lambda(z-\lambda)}\).
The uniform comparisons in \eqref{eq-local-weight} convert the two
integrals of \(|h|^2\) into the corresponding integrals of
\(|V|^2w/w(R)\). Canceling the factor \(1/w(R)\) gives
\begin{equation}\label{s7-weight-annulus}
 \int_{D(\lambda,\delta)}|V(z)|^2w(|z|)\dA(z)
 \le C\int_{D(\lambda,2\delta)\setminus D(\lambda,\delta)}
              |V(z)|^2w(|z|)\dA(z).
\end{equation}
The disks \(D(\lambda,2\delta)\) are pairwise disjoint, and all the
annuli on the right lie in
\(\Omega=\C\setminus\bigcup_{\lambda\in\Lambda_*}D(\lambda,\delta)\).
Summing \eqref{s7-weight-annulus} gives
\[
 \sum_{\lambda\in\Lambda_*}\int_{D(\lambda,\delta)}|V|^2w\,\dA
 \le C\int_\Omega |V|^2w\,\dA.
\]
Adding the integral over \(\Omega\) bounds \(\|V\|_w^2\) by
\(C\int_\Omega|V|^2w\,\dA\).
For \(z\in\Omega\), the distance bound \(d_*(z)\ge\delta\) and
\eqref{eq-F*asymp} give \(w(|z|)\le C|F_*(z)|^2e^{-|z|^2}\).
Substituting into the last integral proves the lower bound in
\eqref{eq-norm-equivalence}. All comparisons apply to the nonnegative
integrals even when one is infinite.
\end{proof}

\begin{proof}[Proof of Theorem~\ref{thm-endpoint}]
By Proposition~\ref{prop-product},
\(\log M_{F_*}(r)/r^2\to1/2\), which gives order two and type \(1/2\).
Equation~\eqref{eq-exponentials} gives \(e^{\lambda z}F_*\in\F_1^2\)
for every \(\lambda\in\C\).
Let \(V(z)=\sum_{k\ge0} a_kz^k\) with \(F_*V\in\F_1^2\).
By Lemma~\ref{lem-norm}, \(\|V\|_w<\infty\). On each circle the Taylor
series converges uniformly, so Parseval's identity in the angular
variable gives
\[
 \frac1{2\pi}\int_0^{2\pi}|V(re^{i\theta})|^2\,d\theta
 =\sum_{k\ge0}|a_k|^2r^{2k}.
\]
We integrate this equality against \(2\pi r w(r)\,dr\).
Tonelli's theorem gives
\begin{equation}\label{s7-moments}
 \|V\|_w^2=\sum_{k\ge0}|a_k|^2\mu_k,\qquad
 \mu_k=2\pi\int_0^\infty r^{2k+1}w(r)\,dr<\infty.
\end{equation}
Each moment \(\mu_k\) is finite by the decay of \(w\).
Apply the same equality to the tail \(V-S_NV\) for
\(S_NV=\sum_{k=0}^Na_kz^k\). The upper bound in Lemma~\ref{lem-norm}
then gives
\begin{equation}\label{eq-tail}
 \|F_*(V-S_NV)\|_1^2
 \le C_0\|V-S_NV\|_w^2
 =C_0\sum_{k>N}|a_k|^2\mu_k\longrightarrow0.
\end{equation}
The limit follows from the summability in \eqref{s7-moments}.
Thus \(\M_1(F_*)\subset[F_*]_1\). The reverse inclusion follows from
the closedness of \(\M_1(F_*)\).
\end{proof}

Under the scaling \eqref{eq-scaling} with \(\alpha=1\), the bound
\eqref{eq-critical-growth} gives
\[
 |F_*(\sqrt\pi\,z)|\le C_A e^{\pi|z|^2/2-A\sqrt\pi\,|z|}
 \qquad(A>0).
\]
Thus \cite[Lemma~4.1]{BB} identifies the closed exponential and
polynomial spans of \(F_*(\sqrt\pi\,\cdot)\) in \(\F_\pi^2\).
The scaling sends polynomials to polynomials and takes the family
of exponential multipliers onto itself. Transforming back, we obtain
\begin{equation}\label{s7-exponential-span}
 \overline{\operatorname{span}\{e^{\lambda z}F_*:\lambda\in\C\}}
 ^{\,\F_1^2}
 =[F_*]_1=\M_1(F_*).
\end{equation}

\noindent \textbf{AI Disclosure.}\enspace
Artificial-intelligence tools were used only for language editing, LaTeX formatting, and limited assistance with local mathematical reasoning. The proof strategy and mathematical development are the authors' own work; all mathematical content was independently written and checked by the authors, who take full responsibility for it.

\end{document}